\documentclass[twoside,12pt]{article}
\usepackage{amssymb,amsmath,bm,euscript}
\usepackage{graphicx,color,caption}

\usepackage[colorlinks=true]{hyperref}

\usepackage{algorithm,algpseudocode}
\usepackage{braket}

\newtheorem{proposition}{Proposition}[section]
\newtheorem{theorem}[proposition]{Theorem}

\newtheorem{corollary}[proposition]{Corollary}

\newcommand{\qed}{\hphantom{.}\hfill $\Box$\medbreak}

\def\F{\mathcal{F}}

\def\x{{\bf x}}
\def\y{{\bf y}}
\def\z{{\bf z}}

\def\uu{{\bf u}}

\def\0{{\bf 0}}

\title{\bf{Three Minimax Ideal Relations of Lie Algebras}\thanks{This research was supported by the Hong Kong
    Research Grant Council (Grant No.  PolyU 15300715, 15301716 and 15300717). }}
\author{ \hspace{1mm} Liqun Qi \thanks{Department of Applied
    Mathematics, The Hong Kong Polytechnic University, Hung Hom,
    Kowloon, Hong Kong; ({\tt liqun.qi@polyu.edu.hk}).} }

\begin{document}
\date{\today}
\maketitle

\begin{abstract}  In this paper, we introduce near perfect ideals and upper bounded ideals, and study them as well as perfect ideals for finite dimensional Lie algebras.  We show that the largest perfect ideal and the largest near perfect ideal of a finite dimensional Lie algebra always exist, and are equal to the smallest ideal of the derived series, and the smallest ideal of the lower central series, respectively.  We call them the perfect radical and the near perfect radical of that Lie algebra, respectively.   A nonzero Lie algebra is solvable if and only if its perfect radical is zero.   The factor algebra of a Lie algebra by its perfect radical is solvable.   A nonzero Lie algebra is nilpotent if and only if its near perfect radical is zero.  The factor algebra of a Lie algebra by its near perfect radical is nilpotent.  We also show that the smallest upper bounded ideal always exists, and is equal to the largest ideal of the upper central series.   For a nilpotent Lie algebra, there is only one upper bounded ideal, i.e., the nilpotent Lie algebra itself.


\end{abstract}


\section{Introduction}

 At the later part of the nineteen century, the Norwegian mathematician Sophus Lie started the research on some continuous transformation groups \cite{LS93}.   Later the theory started by him was called Lie theory, which includes both Lie group theory and Lie algebra theory.   After mathematicians' further works, Lie theory became an important branch of group and algebra theory.   However, it was American physicist Murray Gell-Mann's work which made Lie theory an essential tool for theoretical physicists.   Gell-Mann applied Lie group theory to the study of elementary particles.    In 1961, Gell-Mann applied the $SU(3)$, a Lie group, to his research on tabling elementary particles.   He predicted the Omega minus particle.  This was immediately confirmed by experiments.   Based upon this, in 1964, Gell-Mann further developed the quark theory.   In 1969, he received Nobel prize in physics.   The works of him and other physicists led to a model now called the Standard Model, which is our current best model of the world at the smallest scales.   Gell-Mann's work made theoretical physicists to recognize the importance of group theory, mainly the theory of Lie group and Lie algebra.    Nowadays, group theory in physics is a compulsory course of students of theoretical physics in most universities in the world.

 Lie algebra theory is an algebraic tool for studying Lie groups.   Later it itself became a branch of algebra and found its own application in physics and mechanics \cite{Ia15, LNRW89, RSW97, SW14}.   For mathematicians, a Lie algebra was defined as a vector space with an additional bilinear commutation operation called the Lie bracket operation.   The study on Lie algebra is thus worked more algebraically in a rigid manner \cite{EW06, Ja79}.   Physicists accept this definition, but more prefer to use the structure constants to describe Lie algebras.   This is more workable in computation \cite{Ia15, RWZ88, SW14}.

 Beside Sophus Lie (1842-1899) himself, Friedrich Engel (1861-1941) and Wilhelm Killing (1847-1923) have also made contributions to the early development of Lie algebra.  However, it was \'{E}lie Cartan (1869-1951) who made fundamental contributions to the development of the classical theory of Lie algebra.   He characterized solvable Lie algebras and semisimple Lie algebras by their Killing forms and made a complete classification for semisimple Lie algebras \cite{Ca94}.    Hermann Weyl (1885-1955), Bartel van der Waerden (1903-1996) and Eugene Dynkin (1924-2014) further developed Cartan's work.

 In 1905, Eugenio Elia Levi \cite{Le05} showed that a finite-dimensional Lie algebra is a semidirect sum of a solvable ideal and a semisimple  Lie subalgebra.   This leaves the classification work for finite-dimensional Lie algebras to the classification of solvable Lie algebras and its subclass nilpotent Lie algebras.    There are works on this \cite{GK96, RWZ88, SW14}.  The problem is not totally solved yet.

 The solvable ideal in the Levi decomposition is actually the largest solvable ideal of that Lie algebra, and is called the radical of that Lie algebra.    The largest nilpotent of a Lie algebra also exists, and is called the nilradical of that Lie algebra.

 Not all classes of ideals of a Lie algebra have the largest ideals.  For example, the maximal Abelian ideal of a Lie algebra is not necessarily unique.   The reason for this is that the sum of Abelian ideals is not necessarily an Abelian ideal.

 The radical and the nilradical of a Lie algebra have various important properties and play an important role in classification and identification of Lie algebras.  A nonzero Lie algebra is semisimple if and only its radical is zero.   The factor algebra of a Lie algebra by its radical is semisimple.  The radical of a Lie algebra is the orthogonal complement of the derived algebra of that Lie algebra, with respect to the Killing form of that Lie algebra.

 In the last century, American mathematician Nathan Jacobson (1910-1999) in his book \cite{Ja79}, besides his other contributions to Lie algebra, such as strengthening the Engel theorem \cite{Ja79, Ra87},  studied  radicals and nilradicals, intensively.   Later, more attentions were paid to computing the Levi decomposition and the nilradical, and the classifications of solvable and nilpotent Lie algebras.   In 1988, Rand, Winternitz and Zassenhaus \cite{RWZ88} introduced algorithms for computing the Levi decomposition and the nilradical \cite{SW14}.  The book \cite{GK96} studied more about nilpotent Lie algebras, while the book \cite{SW14} paid more attentions on classification and identification of nilpotent, solvable and Levi decomposable Lie algebras.

We wonder if there are similar structures of a Lie algebra by other kinds of ideals.  In this paper, we introduce near perfect ideals and upper bounded ideals, and study them as well as perfect ideals for finite dimensional Lie algebras.  We show that the largest perfect ideal and the largest near perfect ideal of a finite dimensional Lie algebra always exist, and are equal to the smallest ideal of the derived series, and the smallest ideal of the lower central series, respectively.   Then we introduce upper bounded ideals for finite dimensional Lie algebras.   We show that the smallest upper bounded ideal always exists, and is equal to the largest ideal in the upper central series of that Lie algebra.  These reveal three minimax ideal relations of a Lie algebra.

 In the next section, we study properties of perfect ideals of a finite dimensional Lie algebra.   A Lie algebra is called a perfect Lie algebra if it is equal to its derived algebra.   An ideal of that Lie algebra is called a perfect ideal if it itself is a perfect Lie algebra.   Since zero is a trivial perfect ideal, a Lie algebra always has a perfect ideal.  We show that the sum of two perfect ideals of a Lie algebra is still a perfect ideal of that Lie algebra.   Thus, the largest perfect ideal of a Lie algebra exists.  We show that it is exactly the smallest ideal of the derived series, and call it the perfect radical of that Lie algebra. This reveals a minimax ideal relation of a Lie algebra.   We show that a nonzero Lie algebra is solvable if and only if its perfect radical is zero.  This echoes the known result that a nonzero Lie algebra is semisimple if and only its radical is zero.   We show that if the factor algebra of a Lie algebra by a perfect ideal is also perfect, then that Lie algebra is perfect.   By this, we prove that the factor algebra of a Lie algebra by its perfect radical is solvable.  This echoes the known result that the factor algebra of a Lie algebra by its radical is semisimple.   Some further results on the relations between the perfect radical and the Levi decomposition are presented.   Some further questions on perfect radicals are raised.

In Section 3, we define near perfect ideals of a finite dimensional Lie algebra, and study their properties.   A perfect ideal is a near perfect ideal but not vice versa.    We show that the sum of two near perfect ideals of a Lie algebra is still a near perfect ideal of that Lie algebra.   Thus, the largest near perfect ideal of a Lie algebra exists.  We show that it is exactly the smallest ideal of the lower central series, and call it the near perfect radical of that Lie algebra. This reveals another minimax ideal relation of a Lie algebra.  The near perfect radical contains the perfect radical of a Lie algebra.  We show that a nonzero Lie algebra is nilpotent if and only if its near perfect radical is zero.  This again echoes the known result that a nonzero Lie algebra is semisimple if and only its radical is zero.   We also prove that the factor algebra of a Lie algebra by its near perfect radical is nilpotent.  This again echoes the known result that the factor algebra of a Lie algebra by its radical is semisimple.  Some further questions on near perfect radicals are raised.

We introduce upper bounded ideals for finite dimensional Lie algebras in Section 4.   We show that the smallest upper bounded ideal always exists, and is equal to the largest ideal in the upper central series of that Lie algebra.  We also show that a nilpotent Lie algebra has only one upper bounded ideal, namely the nilpotent Lie algebra itself.

Some further comments are made in Section 5.   In this paper, we only consider finite dimensional Lie algebras.  The related field $\F$ is either the complex field or the real field.

\section{Perfect Ideals of a Lie Algebra}
\label{sec:Lies}

 Suppose that $L$ is an $n$-dimensional Lie algebra defined on $\F$.   The Lie bracket operation on $L$ is denoted by $[ \cdot, \cdot ]$.  Suppose that $I$ and $J$ are two ideals of $L$.  Then $I \cap J$,
 $$I+J := \{ \x + \y : \x \in I, \y \in J \}$$
 and
 $$[I, J] := {\rm Span}\{ [\x, \y] : \x \in I, \y \in J \}$$
 are also ideals of $L$.   The derived algebra of $L$ is defined as $D(L) := [L, L]$.  If $L = D(L)$, then $L$ is called a perfect Lie algebra.   Note that $0$ is a trivial perfect Lie algebra, but there are no one dimensional and two dimensional perfect Lie algebras.  For $\F = \bf C$, up to isomorphism, there is only one perfect Lie algebra in dimension three, namely sl$(2, {\bf C})$, the space of traceless $2 \times 2$ complex matrices \cite{EW06}.

 Let
 $$L^{(0)} = L, \ \ \ L^{(k+1)} = D(L^{(k)}).$$
 Then we have the derived series of $L$: $L^{(0)} \supseteq L^{(1)} \supseteq L^{(2)} \supseteq \cdots$.
 $L^{(k)}$ are ideals of $L$ for all $k$.   If for some $m$, $L^{(m)} = 0$, then $L$ is called solvable.   If $I$ and $J$ are two solvable ideals of $L$, then $I+J$ is also a solvable ideal of $L$.   Since $0$ is a trivial solvable ideal of $L$, the largest solvable ideal of $L$ exists, and is called the radical of $L$ and denoted as $R(L)$.   A nonzero Lie algebra $L$ is called a semisimple Lie algebra if it has no nonzero solvable ideals.  Then, a nonzero Lie algebra is semisimple if and only if its radical is zero.  A semisimple Lie algebra is always a perfect Lie algebra but not vice versa.

 An ideal $I$ of $L$ is called a perfect ideal if it itself is a perfect Lie algebra.

 \begin{proposition}  \label{p2.1}
 Suppose that $L$ is an $n$-dimensional Lie algebra defined on $\F$.   Let $I$ and $J$ be two perfect ideals of $L$.  Then $I+J$ is also a perfect ideal of $L$.
 \end{proposition}
 {\bf Proof} Since $I$ and $J$ are perfect ideals, $I = [I, I] \subseteq [I+J, I+J] = D(I+J)$.   Similarly, $J \subseteq D(I+J)$.  This implies that $I+J \subseteq D(I+J)$.   Hence,
 $I+J = D(I+J)$.   This means that $I+J$ is perfect, hence a perfect ideal of $L$. \qed

 Since zero is a trivial perfect ideal of $L$, the largest perfect ideal of $L$ exists.   We call it the perfect radical of $L$ and denote it as $P(L)$.

  \begin{proposition}  \label{p2.2}
 A nonzero Lie algebra is solvable if and only if its perfect radical is zero.
 \end{proposition}
 {\bf Proof} If $L$ is not solvable, then either it itself is perfect, or there is an $m$ such that
$$L^{(m-1)} \supsetneq L^{(m)} = L^{(m+1)} \not = 0.$$
In the first case, $L$ itself is a nonzero perfect ideal.   In the second case, $L^{(m)}$ is a nonzero perfect ideal.
On the other hand, suppose that $I$ is a nonzero perfect ideal of $L$.  Then
$$I = D(I) \subseteq D(L) = L^{(1)},$$
and by induction we have
$$I = D(I) \subseteq D(L^{(k)}) = L^{(k+1)}$$
for all $k$.   Hence, $L$ cannot be solvable. \qed

We now reveal a minimax ideal relation of $L$.

\begin{proposition} \label{p2.3}
The perfect radical of $L$ is equal to the smallest ideal of the derived series of $L$.
\end{proposition}
{\bf Proof}  We may discuss this in three cases.

The first case is that $L$ is solvable.  In this case the smallest ideal of the derived series is zero.  On the other hand, in this case, $P(L)=0$ by Proposition \ref{p2.2}.   Hence, $P(L)$ is equal to the smallest ideal of the derived series in this case.

The second case is that $L$ is perfect.  In this case the smallest ideal of the derived series is $L$ itself.   On the other hand, in this case, $P(L)=L$ by the definition of perfect Lie algebras.   Hence, $P(L)$ is equal to the smallest ideal of the derived series in this case too.

The third case is that $L$ is neither solvable nor perfect.  Then, consider
$L^{(m)}$, where $L^{(m)}$ satisfies
\begin{equation} \label{e2.1}
L^{(m-1)} \supsetneq L^{(m)} = L^{(m+1)} \not = 0
\end{equation}
in the derived series.   Then $L^{(m)}$ is the smallest ideal of the derived series.
By (\ref{e2.1}), $L^{(m)}$ is a perfect ideal of $L$.  Hence,
$$L^{(m)} \subseteq P(L).$$
Since $P(L) \subseteq L$,
$$P(L) = D(P(L)) \subseteq D(L) = L^{(1)}.$$
By induction, we have
$$P(L) = D(P(L)) \subseteq D(L^{(k)}) = L^{(k+1)}$$
for all $k$.  Thus,
$$P(L) \subseteq L^{(m)}.$$
This shows that
$$P(L) = L^{(m)},$$
i.e., $P(L)$ is equal to the smallest ideal of the derived series in all the cases.
\qed

\medskip

Suppose that $I$ is an ideal of $L$.  Then the quotient vector space $L/I = \{ \x + I : \x \in L \}$ is a Lie algebra with a Lie bracket on $L/I$ defined by
$$[\x + I, \y + I] := [\x, \y] + I, \ \ \forall \x, \y \in L,$$
and is called the quotient or factor algebra of $L$ by $I$.

\begin{proposition}  \label{p2.4}
 Suppose that $L$ is an $n$-dimensional Lie algebra defined on $\F$, and has a perfect ideal $I$.  If furthermore $L/I$ is perfect, then $L$ is a perfect Lie algebra.
\end{proposition}
{\bf Proof} Let $\x \in L$.   Since $L/I$ is perfect, there are finitely many $\y_i, \z_i \in L$, $i = 1, \cdots, s$, such that
$$\x + I = \sum_{i=1}^s [\y_i + I, \z_i + I].$$
This implies that
$$\x = \sum_{i=1}^s [\y_i, \z_i] + \uu$$
for some $\uu \in I$.   Since $I$ is perfect,
$$\uu \in I = [I, I] \subseteq [L, L] = D(L).$$
On the other hand, $[\y, \z] \in D(L)$.   Thus, $\x \in D(L)$.   As $\x$ is arbitrary in $L$, we have $L \subseteq D(L)$, i.e., $L = D(L)$ and $L$ is perfect.
\qed

\medskip

Recall that $L/R(L)$ is always semisimple.   We have a similar result for the perfect radical $P(L)$.

\begin{proposition}  \label{p2.5}
 Suppose that $L$ is an $n$-dimensional Lie algebra defined on $\F$.
 Then $L/P(L)$, the factor algebra of $L$ by its perfect radical $P(L)$,  is solvable.
\end{proposition}
{\bf Proof} Let $H$ be a perfect ideal of $L/P(L)$.   By the ideal correspondence, there is an ideal $J$ of $L$ containing $P(L)$ such that $H = J/P(L)$.  By definition, $P(L)$ is perfect, and $J/P(L) = H$ is perfect by hypothesis.  Therefore Proposition \ref{p2.4} implies that $J$ is perfect.  But then $J$ is contained in $P(L)$; that is, $H = 0$. By Proposition \ref{p2.3}, $L/P(L)$ is solvable. \qed

\medskip
We have the following question.

{\bf Question 1} Suppose that $L$ is an $n$-dimensional Lie algebra defined on $\F$.  Is $L$ always decomposable to the semidirect sum of $P(L)$ and a solvable subalgebra of $L$?

If the answer to this question is ``yes'', then we have a decomposition of a general finite dimensional Lie algebra $L$, other than the Levi decomposition.

\medskip

According to the Levi theorem, any finite dimensional Lie algebra $L$ can be decomposed into the semidirect sum of its radical $R(L)$ and a semisimple Lie subalgebra $S$ of $L$.
The semisimple Lie subalgebra $S$ is called a Levi factor of $L$.   The Levi factor is not unique.  However, any two Levi factors are isomorphic in the sense of the Mal'cev theorem \cite{Ja79, SW14}.

Let
 $$L^0 = L, \ \ \ L^{k+1} = [L, L^k].$$
 Then we have the lower central series of $L$: $L^0 \supseteq L^1 \supseteq L^2 \supseteq \cdots$. $L^k$ are ideals of $L$ for all $k$.   If for some $m$, $L^m = 0$, then $L$ is called nilpotent.  A nilpotent Lie algebra is always a solvable Lie algebra but not vice versa. A Lie algebra $L$ is solvable if and only if $D(L)$ is nilpotent.   If $I$ and $J$ are two nilpotent ideals of $L$, then $I+J$ is also a nilpotent ideal of $L$.   Since $0$ is a trivial nilpotent ideal of $L$, the largest nilpotent ideal of $L$ exists, and is called the nilradical of $L$ and denoted as $N(L)$.   Since a nilpotent
Lie algebra is solvable, we always have $N(L) \subseteq R(L)$.   By Theorem 13 of \cite{Ja79}, we have $[L, R(L)] \subseteq N(L)$.   Hence, we always have
\begin{equation} \label{e2.2}
[L, R(L)] \subseteq N(L) \subseteq R(L).
\end{equation}

\begin{theorem} \label{t2.6}
Suppose that $L$ is an $n$-dimensional Lie algebra defined on $\F$, $P(L)$ is its perfect radical, $S$ is a Levi factor of $L$.  Then we have the following conclusions.

(a) $S \subseteq P(L)$.

(b) $S$ is also a Levi factor of $P(L)$.

(c) $R(P(L)) = R(L) \cap P(L) \subseteq N(P(L))$.   Thus $R(P(L)) = N(P(L))$ and is nilpotent.
\end{theorem}
{\bf Proof} If $L$ is solvable, then $P(L) = 0$, $R(L) = L$ and $S = 0$, then (a), (b) and (c) hold.   If $L$ is perfect, then $P(L) = L$, (a), (b)  and $R(P(L)) = R(L) \cap P(L)$
 also hold.  Now assume that $L$ is neither solvable nor perfect.  Then this is the third case of the proof of Proposition \ref{p2.3}.   Let $L^{(k)}$ be as in the third case of Proposition \ref{p2.3} for $k = 0, \cdots, m$.  Since $S$ is a Levi factor, it is semisimple, and hence perfect.   We have $[S, S] = S$.  Since $S \subseteq L$, we have
$$S = [S, S] \subseteq [L, L] = L^{(1)}.$$
By induction, we have
$$S = [S, S] \subseteq [L^{(k)}, L^{(k)}] = L^{(k+1)}$$
for all $k$.  Hence,
$$S \subseteq L^{(m)} = P(L).$$
Then $S$ must be also a Levi factor of $P(L)$.  This implies (b) and $R(P(L)) = R(L) \cap P(L)$ hold in this case.

We now show that in both the second and the third cases, we have
$$R(P(L)) \subseteq N(P(L)).$$
We have
$$P(L) = S \oplus R(P(L)).$$
Then
\begin{eqnarray*}
P(L) & = & [P(L), P(L)] \\
& \subseteq & [S +R(P(L)), S +R(P(L))] \\
& \subseteq & [S, S]  + [P(L), R(P(L))] \\
& = & S  +  [P(L), R(P(L))] \\
& \subseteq & S  + N(P(L)),
\end{eqnarray*}
where we use $[P(L), R(P(L))] \subseteq  N(P(L))$ by (\ref{e2.2}).
Then we have
\begin{eqnarray*}
R(P(L)) & \subseteq & R(L) \cap [S  + N(P(L))] \\
& = & [S +R(P(L)), S +R(P(L))] \\
& \subseteq & R(L) \cap N(P(L)) \\
& \subseteq & N(P(L)).
\end{eqnarray*}
This implies that $R(P(L)) = N(P(L))$ and is nilpotent.

Hence,  (a), (b) and (c) hold in all the three cases.  \qed

\begin{corollary}
If a finite dimensional Lie algebra $L$ is perfect, then its radical and nilradical are the same, i.e., $R(L)=N(L)$.
\end{corollary}
{\bf Proof} Since $L = P(L)$, this follows from Theorem \ref{t2.6} (c).

\medskip

{\bf Question 2} Is $P(L)$ the union of all Levi factors of $L$?

\medskip

Let gl$(L)$ be the space of all the linear transformations of $L$.   The adjoint representation ad $L$ is a linear map of $L$ into gl$(L)$:
$${\rm ad} : L \to {\rm gl}(L) : \x \to {\rm ad}(\x)$$
defined for any $\x, \y \in L$ via
$${\rm ad}(\x) \y = [\x, \y].$$
The Killing form $K$ of $L$ is a symmetric bilinear form on $L$ defined by
$$K(\x, \y) = {\rm tr}({\rm ad}\left(\x) \cdot {\rm ad}(\y)\right),$$
where tr denotes the trace of the square matrix corresponding the linear transformation in the brackets.
The Killing form plays a fundamental role in Cartan's criteria for solvable Lie algebras and semisimple Lie algebras.

Suppose $I$ is an ideal of $L$.  Then the orthogonal complement of $I$ with respect to the Killing form
$$I^\perp = \{ \x \in L : K(\x, \y) = 0, \forall \y \in I \}$$
is again an ideal of $L$.
By Theorem 5 of \cite{Ja79},  the radical of a Lie algebra is the orthogonal complement of the derived algebra of that Lie algebra with respect to the Killing form:
$$R(L) = D(L)^\perp.$$
This gives a computable formula of the radical $R(L)$ \cite[(6.6)]{SW14}.

We have the following proposition.

\begin{proposition}  \label{p2.8}
 Suppose that $L$ is an $n$-dimensional Lie algebra defined on $\F$.  Then
 $$N(L) \subseteq L^\perp \equiv \{ \x \in L : K(\x, \y) = 0, \forall \y \in L \}.$$
 However, the converse is not true, i.e., $L^\perp$ may not be nilpotent.
\end{proposition}
{\bf Proof} Denote $H= N(L)$.   Since $H$ is an ideal of $L$, $H^k$ is also an ideal of $L$ for all positive integer $k$.  Let $\x \in H$ and $\y \in L$.  Then
$${\rm ad}\x \cdot {\rm ad} \y (L) \subseteq H.$$
Then for $k = 0, 1, 2, \cdots$, we have
$${\rm ad}\x \cdot {\rm ad}\y (H^k) \subseteq H^{k+1}.$$
Since $H$ is nilpotent, $K(\x, \y) = 0$.  This shows that $N(L)\subseteq L^\perp$.

On the other hand, though the Killing form of a nilpotent Lie algebra is zero, the converse is not true.   Such a Lie algebra is solvable but not nilpotent.   Let $L$ be such a Lie algebra.  Then $L^\perp = L$, which is solvable but not nilpotent. \qed
\medskip

{\bf Question 3} Is the Killing form of $L^\perp = L$ always zero?

\medskip

{\bf Question 4}  What is $P(L)^\perp$?

\section{Near Perfect Ideals}

Suppose that $L$ is an $n$-dimensional Lie algebra defined on $\F$.   An ideal $I$ of $L$ is called a near perfect ideal of $L$ if $[L, I] = I$.  Then a perfect ideal is always a near perfect ideal but not vice versa.  Let $\F = \bf C$ and $n = 3$.   Let $\{ \x, \y, \z \}$ be a basis of $L$.   Suppose that
$$[\x, \y] = \0, \ \ [\z, \x] = \x, \ \ [\z, \y] = \x + \y.$$
This is a solvable Lie algebra, denoted as $s_{3, 2}$ in \cite[Page 226]{SW14}. Let $I$ be the subspace spanned by $\x$ and $\y$.  Then $I$ is not a perfect ideal of $L$.  Actually, as we said before, there are no two dimensional perfect Lie algebras.  Hence, there are no two dimensional perfect ideals.   On the other hand, by definition, we see that $I$ is a near perfect ideal of $L$.

\begin{proposition}  \label{p3.1}
 Suppose that $L$ is an $n$-dimensional Lie algebra defined on $\F$.   Let $I$ and $J$ be two near perfect ideals of $L$.  Then $I+J$ is also a near perfect ideal of $L$.
 \end{proposition}
 {\bf Proof} Since $I$ and $J$ are near perfect ideals,
 $$[L, I+J] = [L, I] + [L, J]= I+J.$$
 Hence, $I+J$ is also a near perfect ideal of $L$. \qed

 Since zero is a trivial near perfect ideal of $L$, the largest near perfect ideal of $L$ exists.   We call it the near perfect radical of $L$ and denote it as $NP(L)$.  Since a perfect ideal is always a near perfect ideal, we always have
 $$P(L) \subseteq NP(L).$$

  \begin{proposition}  \label{p3.2}
 A nonzero Lie algebra is nilpotent if and only if its near perfect radical is zero.
 \end{proposition}
 {\bf Proof} If $L$ is not nilpotent, then either it itself is perfect, or there is an $m$ such that
$$L^{m-1} \supsetneq L^m = L^{m+1} \not = 0.$$
In the first case, $L$ itself is a nonzero near perfect ideal.   In the second case, $L^m$ is a nonzero near perfect ideal.
On the other hand, suppose that $I$ is a nonzero near perfect ideal of $L$.  Then
$$I = [L, I] \subseteq [L, L] = L^1,$$
and by induction we have
$$I = [L, I] \subseteq [L, L^k] = L^{k+1}$$
for all $k$.   Hence, $L$ cannot be nilpotent. \qed

We now reveal another minimax ideal relation of $L$.

\begin{proposition} \label{p3.3}
The near perfect radical of $L$ is equal to the smallest ideal of the lower central series of $L$.
\end{proposition}
{\bf Proof}  We may discuss this in three cases.

The first case is that $L$ is nilpotent.  In this case the smallest ideal of the lower central series is zero.  On the other hand, in this case, $NP(L)=0$ by Proposition \ref{p3.2}.   Hence, $NP(L)$ is equal to the smallest ideal of the lower central series in this case.

The second case is that $L$ is perfect.  In this case the smallest ideal of the lower central series is $L$ itself.   On the other hand, in this case, $NP(L)=L$ by the definition of perfect Lie algebras.   Hence, $NP(L)$ is equal to the smallest ideal of the lower central series in this case too.

The third case is that $L$ is neither nilpotent nor perfect.  Then, consider
$L^m$, where $L^m$ satisfies
\begin{equation} \label{e3.3}
L^{m-1} \supsetneq L^m = L^{m+1} \not = 0
\end{equation}
in the lower central series.   Then $L^m$ is the smallest ideal of the lower central series.
By (\ref{e3.3}), $L^m$ is a near perfect ideal of $L$.  Hence,
$$L^m \subseteq NP(L).$$
Since $NP(L) \subseteq L$,
$$NP(L) = [L, NP(L)] \subseteq [L, L] = L^1.$$
By induction, we have
$$NP(L) = [L, NP(L)] \subseteq [L, L^k] = L^{k+1}$$
for all $k$.  Thus,
$$NP(L) \subseteq L^m.$$
This shows that
$$NP(L) = L^m,$$
i.e., $NP(L)$ is equal to the smallest ideal of the lower central series in all the cases.
\qed

We also have the following proposition.

\begin{proposition}  \label{p3.4}
 Suppose that $L$ is an $n$-dimensional Lie algebra defined on $\F$, and has a near perfect ideal $I$ and an ideal $J$ such that $I \subseteq J$.  If furthermore $J/I$ is a near perfect ideal of $L/I$, then $J$ is also a near perfect ideal of $L$.
\end{proposition}
{\bf Proof} Let $\x \in J$.   Since $J/I$ is a near perfect ideal of $L/I$, there are finitely many $\y_i \in L$ and $\z_i \in J$, $i = 1, \cdots, s$, such that
$$\x + I = \sum_{i=1}^s [\y_i + I, \z_i + I].$$
This implies that
$$\x = \sum_{i=1}^s [\y_i, \z_i] + \uu$$
for some $\uu \in I$.   Since $I$ is a near perfect ideal of $L$,
$$\uu \in I = [L, I] \subseteq [L, J].$$
On the other hand, $[\y, \z] \in [L, J]$.   Thus, $\x \in [L, J]$.   As $\x$ is arbitrary in $J$, we have $J \subseteq [L, J]$, i.e., $J = [L, J]$ and $J$ is also a near perfect ideal of $L$.
\qed

Then we have the following proposition similar to Proposition \ref{p2.5}.

\begin{proposition}  \label{p3.5}
 Suppose that $L$ is an $n$-dimensional Lie algebra defined on $\F$.
 Then $L/NP(L)$, the factor algebra of $L$ by its near perfect radical $NP(L)$,  is nilpotent.
\end{proposition}
{\bf Proof} Let $H$ be a near perfect ideal of $L/NP(L)$.   By the ideal correspondence, there is an ideal $J$ of $L$ containing $NP(L)$ such that $H = J/NP(L)$.  By definition, $NP(L)$ is a near perfect ideal of $L$, and $J/NP(L) = H$ is a near perfect ideal of $L/NP(L)$ by hypothesis.  Therefore Proposition \ref{p3.4} implies that $J$ is a near perfect ideal of $L$.  But then $J$ is contained in $NP(L)$; that is, $H = 0$. By Proposition \ref{p3.3}, $L/NP(L)$ is nilpotent. \qed

\medskip
Suppose that $L$ is an $n$-dimensional Lie algebra defined on $\F$.  There are some further questions on its near perfect radical.

\medskip

{\bf Question 5} Is $L$ always decomposable to the semidirect sum of $NP(L)$ and a nilpotent subalgebra of $L$?

\medskip

{\bf Question 6}  What is $NP(L)^\perp$?

\section{Upper Bounded Ideals of a Lie Algebra}

Suppose that $L$ is an $n$-dimensional Lie algebra defined on $\F$.   Let $I$ be an ideal of $L$.   Define
$$U(I) = \{ \x \in L : [\x, L] \subseteq I \}.$$
Then $U(I)$ is also an ideal of $L$.   We have $I \subseteq U(I) \subseteq L$.    We call $U(I)$ the upper extension of $I$.  If $U(I) = I$, then we call $I$ an upper bounded ideal of $L$.   Then $L$ itself is an upper bounded ideal of $L$.

\begin{proposition}
Suppose that $L$ is an $n$-dimensional Lie algebra defined on $\F$, $I$ and $J$ are two upper bounded ideals of $L$.  Then $H=I\cap J$ is also an upper bounded ideal of $L$.
\end{proposition}
{\bf Proof} We know that $H$ is an ideal of $L$.  Let $\x \in L \setminus H$.   Then either $\x \in L \setminus I$ or $\x \in L \setminus J$, or both.  If $\x \in L \setminus I$, then $[\x, L] \not \subset I$ since $I$ is an upper bounded ideal.  If $\x \in L \setminus J$, then $[\x, L] \not \subset J$ since $J$ is an upper bounded ideal.   Then either $[\x, L] \not \subset I$ or $[\x, L] \not \subset J$.  This implies that
$$[\x, L] \not \subset I \cap J = H.$$
Hence, $H$ is an upper bounded ideal of $L$.
\qed

Therefore£¬the smallest upper bounded ideal of $L$ always exist.

Denote $U_1(0):=U(0)$ and $U_{k+1}(0) := U(U_k(0))$ for all positive integer $k$.   We have $U_1(0) = Z(L)$.     Then either $0$ is an upper bounded ideal of $L$, or we have a positive integer $m$ such that
\begin{equation} \label{e4.4}
0 \subsetneq U_1(0)\subsetneq \cdots \subsetneq U_m(0) \equiv H,
\end{equation}
where $U_m(0) = H$ is an upper bounded ideal of $L$.   The ideal series (\ref{e4.4}) is the upper central series of $L$. It forms the three characteristic series of $L$, with the derived series of $L$, and the lower central series of $L$ together.

We now reveal the third minimax ideal relation of $L$.

\begin{proposition}
Suppose that $L$ is an $n$-dimensional Lie algebra defined on $\F$.   Then the smallest
upper bounded ideal of $L$ is $U_m(0)$ in (\ref{e4.4}), the largest ideal in the upper central series of $L$.
\end{proposition}
{\bf Proof} Let $I$ be an upper bounded ideal of $L$.   Then $0 \subseteq I$.  We have
$$U_1(0) = U(0) \subseteq U(I) = I.$$
If $U_k(0) \subseteq I$, then
$$U_{k+1}(0) = U(U_k(0)) \subseteq U(I) = I.$$
Hence, by induction, we have
$$U_m(0) \subseteq I.$$
Since $U_m(0)$ is also an upper bounded ideal of $L$, $U_m(0)$ is the smallest upper bounded ideal of $L$.
\qed

We now show that a nilpotent Lie algebra $L$ has only one upper bounded ideal, i.e., $L$ itself.

\begin{theorem} \label{t4.3}
Suppose that $L$ is an $n$-dimensional nilpotent Lie algebra on $\F$, where $n$ is a positive integer.   Then $L$ has only one upper bounded ideal, i.e., $L$ itself.
\end{theorem}
{\bf Proof} If $n=1$, then $L$ is an Abelian Lie algebra.  It has only two ideals, $0$ and $L$.  We have $U(0) = L$.  Thus there is only one upper bounded ideal, namely $L$ itself.

We now assume that $n \ge 2$.  Let $I$ be an ideal of $L$, $I \not = L$. Suppose that $I$ is an ideal with dimension $n-1$.   Let $\x \in L \setminus I$.   Then
$$[\x, L] \subseteq {\rm Span}\{ [\x, \x] + [\x, I] \} \subseteq [\x, I] \subset I.$$
Hence, $\x \in U(L)$.   This implies that $U(I) = L$. i.e., $I$ is not an upper bounded ideal of $L$.  We now assume that  dim$(I) \le n-2$.   Consider the factor algebra $L/I$.   Then $L/I$ is a nilpotent Lie algebra since $L$ is a nilpotent Lie algebra.  Then the center $Z(L/I)$ is not trivial.   Suppose that $\x+I \in Z(L/I)$ and $\x \not \in I$, i.e., $\x+I$ is not the zero element of $L/I$. Then $[\x+I, \y+I] = I$ for any $\y \in L$. However, $[\x+I, \y+I] = [\x, \y] +I$.   This implies $[\x, \y] + I = I$, i.e., $[\x, \y] \in I$ for any $\y \in L$.   This implies $\x \in U(I)$.   Hence, $I \not = U(I)$.    Therefore, $I$ is not an upper bounded ideal of $L$.  Thus, $L$ has only one upper bounded ideal, namely $L$ itself.
\qed

\begin{corollary} \label{c4.4}
 Suppose that $L$ is a nilpotent Lie algebra.   Then the upper central series of $L$ has the form:
  \begin{equation} \label{e4.5}
0 \subsetneq U_1(0) \equiv Z(L) \subsetneq \cdots \subsetneq U_m(0) = L.
\end{equation}
\end{corollary}

\section{Further Comments}

In this paper, we introduce near perfect ideals and upper bounded ideals, and study them as well as perfect ideals for finite dimensional Lie algebras.  We show that the largest perfect ideal and the largest near perfect ideal of a finite dimensional Lie algebra always exist, and are equal to the smallest ideal of the derived series, and the smallest ideal of the lower central series, respectively. We also show that the smallest upper bounded ideal always exists, and is equal to the largest ideal of the upper central series.   These reveal three minimax ideal relations of a Lie algebra.    Minimax relations are some essential properties in many branches of mathematics.  It is worth further exploring such relations in Lie algebra.

We call the largest perfect ideal and the largest near perfect ideal of a Lie algebra, the perfect radical and the near perfect radical of that Lie algebra, respectively. We discover that the relation between solvable Lie algebras and perfect radicals, and the relation between nilpotent Lie algebras and near perfect radicals, are very similar to the relation between semisimple Lie algebras and radicals.     We thus think that there may be some structural essence of Lie algebras behind this similarity, and plan to explore the properties of
perfect radicals and near perfect radicals further.

\bigskip

{\bf Acknowledgment}   The author is thankful to Professors Chengming Bai, Huixiang Chen, Shenglong Hu, Libor \u{S}nobl, Yisheng Song and Dr. Yi Xu for their comments.

\end{document}